\documentclass{compositio}

\usepackage{amsmath} 

\newtheorem{defin}{D\'ef\mbox{}inition}[section]
\newtheorem{theo}[defin]{Th\'eor\`eme}
\newtheorem{prop}[defin]{Proposition}
\newtheorem{lem}[defin]{Lemme}
\newtheorem{cor}[defin]{Corollaire}
\newtheorem{rem}[defin]{Remarque}

\newtheorem*{defin*}{D\'ef\mbox{}inition}
\newtheorem*{theo*}{Th\'eor\`eme}
\newtheorem*{prop*}{Proposition}
\newtheorem*{lem*}{Lemme}
\newtheorem*{cor*}{Corollaire}
\newtheorem*{rem*}{Remarque}
\newtheorem*{aff*}{Affirmation}
\newtheorem*{affs*}{Affirmations}
\newtheorem*{nota*}{Notation}

\newtheorem*{theobis*}{Th\'eor\`eme 2 bis}

\newtheorem{theoi}{Th\'eor\`eme}

\newcommand{\dps}{\displaystyle}

\newcommand{\gauche}{\begin{flushleft}\end{flushleft}}

\newcommand{\R}{\mathbb{R}}
\newcommand{\C}{\mathbb{C}}
\newcommand{\N}{\mathbb{N}}
\newcommand{\Z}{\mathbb{Z}}
\newcommand{\Q}{\mathbb{Q}}

\newcommand{\An} {\mathbf{A}_n}

\newcommand{\D}{\mathcal{D}}
\newcommand{\Dn}{\mathcal{D}_n}
\newcommand{\Dnp}{\mathcal{D}_{n+p}}

\newcommand{\Dnz}{\mathcal{D}_n \langle z \rangle}
\newcommand{\Dnpz} {\mathcal{D}_{n+p}\langle z \rangle}
\newcommand{\Dm}{\mathcal{D}_m}
\newcommand{\Dmz}{\mathcal{D}_m\langle z \rangle}
\newcommand{\Dnd}{\mathcal{D}_{n+2}}

\newcommand{\dxi}{\partial _{x_i}}

\newcommand{\dx}[1]{\partial _{x_{#1}}}
\newcommand{\dt}[1]{\partial _{t_{#1}}}
\newcommand{\ddx}{\partial _x}
\newcommand{\ddt}{\partial _t}

\newcommand{\dxsur}[2]{\frac{\partial {#1}}{\partial x_{#2}}}

\newcommand{\B}{\mathcal{B}}

\newcommand{\ord}{\text{ord}}
\newcommand{\gr}{\text{gr}}
\renewcommand{\mp}{\text{mp}} 

\newcommand{\DN}{\mathcal{N}} 
\newcommand{\Exp}{\text{Exp}}

\newcommand{\Vb}{\overline{V}}
\newcommand{\U}{\mathcal{U}}
\newcommand{\UV}{\mathcal{U}_V}
\newcommand{\E}{\mathcal{E}}
\newcommand{\EV}{\mathcal{E}_V}
\renewcommand{\L}{\mathcal{L}}

\begin{document}

\title[Existence de polyn\^omes de Bernstein-Sato~: une d\'emonstration
  constructive]{D\'emonstration constructive de l'existence de
  polyn\^omes de Bernstein-Sato pour plusieurs fonctions analytiques}

\author{Rouchdi BAHLOUL}
\email{rouchdi.bahloul@univ-angers.fr}
\address{D\'epartement de Math\'ematiques\\U.M.R. 6093\\Universit\'e
  d'Angers\\2 bd Lavoisier, 49045 Angers cedex 01, France}




\begin{abstract}
In 1987, C.~Sabbah proved the existence of Bernstein-Sato polynomials
associated with several analytic functions. The purpose of this
article is to give a more elementary and constructive proof of the
result of C.~Sabbah based on the notion of the analytic Gr\"obner fan
of a $\D$-module.
\end{abstract}

\maketitle




\section*{Introduction et \'enonc\'e des r\'esultats principaux}

Fixons $n\ge 1$ et $p\ge 1$ deux entiers et $v\in \N^p$. Soient
$x=(x_1,\ldots,x_n)$ et $s=(s_1,\ldots,s_p)$ deux syst\`emes de
variables. On se donne $f_1,\ldots, f_p \in \C\{x\}=\C\{x_1,\ldots,
x_n\}$. Notons $\Dn$ l'anneau des op\'erateurs diff\'erentiels \`a
coefficients dans $\C\{x\}$.
Pour $b(s) \in \C[s]=\C[s_1,\ldots,s_p]$, consid\'erons l'identit\'e
suivante~: 
\[(\star) \qquad b(s) f^s \in \D_n[s] f^{s+v},\]
o\`u $f^{s+v}=f_1^{s_1+v_1} \cdots f_p^{s_p+v_p}$. Un polyn\^ome
$b(s)$ v\'erifiant une telle identit\'e est appel\'e polyn\^ome de
Bernstein-Sato (associ\'e \`a $f=(f_1,\ldots,f_p)$). L'ensemble de ces
polyn\^omes forment un id\'eal appel\'e id\'eal de Bernstein-Sato
(associ\'e \`a $f$) et qu'on note $\B^v(f)$.\\
Rappelons que c'est I.N. Bernstein \cite{bernstein} qui, dans le cas
$p=1$ et o\`u $f$  est polynomiale, a montr\'e que l'id\'eal $\B^v(f)$
est non nul (dans ce cas, il faut, dans $(\star)$, remplacer $\D_n$
par l'alg\`ebre de Weyl $\An(\C)$, i.e. l'alg\`ebre des op\'erateurs
diff\'erentiels \`a coefficients polynomiaux). Dans le cas o\`u $f$
est analytique et toujours pour $p=1$, la non nullit\'e de $\B^v(f)$
revient \`a J.E.~Bj\"ork \cite{bjork} avec des m\'ethodes similaires
\`a celles employ\'ees dans \cite{bernstein}. Dans ce m\^eme cas,
citons M.~Kashiwara \cite{kashiwara} qui publia une autre preuve et
d\'emontra en plus que le g\'en\'erateur unitaire de l'id\'eal de
Bernstein-Sato est \`a racines rationnelles.
Pour $p\ge 2$, la preuve dans le cas polynomial est une
g\'en\'eralisation facile de celle de I.N.~Bernstein, que l'on peut
trouver dans \cite{lichtin}. Dans le cas analytique avec $p\ge 2$, la
non nullit\'e de $\B^v(f)$ a \'et\'e d\'emontr\'ee par C.~Sabbah
(\cite{sabbah1} et \cite{sabbah2}). Citons la contribution de A.~Gyoja
\cite{gyoja} qui a repris la preuve de C.~Sabbah en montrant de plus
que $\B^v(f)$ contient un \'el\'ement rationnel non nul.

L'objet du pr\'esent article est une mise au point de la preuve
donn\'ee par C. Sabbah.
Plus pr\'ecis\'ement, on peut d\'ecomposer la preuve de C.~Sabbah en
deux grandes \'etapes~: la premi\`ere utilise des arguments similaires
\`a ceux employ\'es par M.~Kashiwara dans le cas $p=1$, la deuxi\`eme
consiste essentiellement en un r\'esultat de finitude qui permet par
la suite de se ramener au r\'esultat de la premi\`ere \'etape. 
La seconde \'etape de la preuve de C.~Sabbah s'appuie sur un
\'eventail dit adapt\'e (\cite{sabbah1} prop. 2.2.1 et th. A.1.1) dont
l'existence m\'eriterait une mise au point technique
suppl\'ementaire (voir le commentaire qui suit le th.~S1). Aussi, nous
proposons dans cet article un \'enonc\'e et une d\'emonstation plus
\'el\'ementaires et plus constructifs de la seconde \'etape, qui
\'evitent la notion d\'elicate d'\'eventail adapt\'e.

Afin de motiver les r\'esultats du pr\'esent article, nous rappelons,
sans entrer dans tous les d\'etails, la preuve de C.~Sabbah. Dans
ce rappel, nous mettrons en \'evidence le r\'esultat de C.~Sabbah pour
lequel nous allons donner un \'enonc\'e et une preuve plus
constructifs. Signalons que la plupart des notions ou notations
introduites ci-dessous qui nous seront utiles dans la suite seront
d\'etaill\'ees dans les sections suivantes.

On note $\Dnp$ l'anneau des op\'erateurs diff\'erentiels \`a
coefficients dans $\C\{x,t\}=\C\{x_1,\ldots,x_n,\\ t_1,\ldots,t_p\}$. En
suivant la m\'ethode de B.~Malgrange \cite{malgrange}, on fait agir
$\Dnp$ sur $\C\{x\}[\frac{1}{f_1\cdots f_p}, s] f^s$. On note $I$
l'id\'eal (\`a gauche\footnote{Dans tout le texte, id\'eal signifie
  id\'eal \`a gauche}) annulateur de $f^s$ dans $\Dnp$ et $M$ le
quotient $M=\Dnp/I$.

Pour chaque $j=1,\ldots,p$, notons $V_j(\Dnp)$ la $V$-filtration de
Kashiwara-Malgrange associ\'ee \`a la variable $t_j$ et notons
$V=(V_1,\ldots,V_p)$ la (multi)filtration de $\Dnp$ index\'ee par
$\Z^p$~: pour $w\in \Z^p$,
\[V_w(\Dnp)=\bigcap_{j=1}^p \{V_j\}_{w_j}(\Dnp).\]
Elle induit une filtration $V(M)$ sur $M$ o\`u pour tout $w\in \Z^p$,
$V_w(M)$ est l'image de $V_w(\Dnp)$ par la projection $\Dnp \to
M=\Dnp/I$.

Pour $j=1,\ldots,p$, identifions la filtration $V_j$ avec la forme
lin\'eaire sur $\N^{2n+2p}$ donn\'ee par $V_j(\alpha,\mu,
\beta,\nu)=\nu_j -\mu_j$ (o\`u $\alpha,\beta \in \N^n$, $\mu,\nu \in
\N^p$ et $\alpha, \mu, \beta, \nu$ correspondent respectivement \`a
$x, t, \ddx, \ddt$).
Notons $\UV=\sum_{j=1}^p \R_{\ge 0} V_j$. On identifie $\UV$ \`a
$(\R_{\ge 0})^p$. Chaque $L$ de $\UV \cap \N^p$ (i.e. $L$ \`a
coefficients entiers) donne lieu \`a une filtration
naturelle $V^L$ de $\Dnp$ et de $M$ index\'ee par $\Z$ donn\'ee par~:
\[V^L_k(M) = \sum_{\{w\in \Z^p;\, L(w)\le k\} } V_w(M),\]
o\`u $L(w)=l_1 w_1+\cdots+l_p w_p$ si $L=(l_1,\ldots,l_p) \in \N^p$.\\
\\
Voici maintenant les deux \'etapes de la preuve donn\'ee par C.~Sabbah.
\begin{description}
\item[Etape 1]
\begin{theo*}(\cite{sabbah1} th. 3.1.1, voir aussi \cite{gyoja} 2.9 et
  2.10)
Pour tout $L \in \UV \cap \N^p$, il existe un polyn\^ome $b \in
\C[\lambda]$ d'une variable tel que pour tout $k \in \Z$, on ait
\[b\big(L(-\dt{1} t_1, \ldots, -\dt{p}t_p)-k\big) V^L_k(M) \subset
V^L_{k-1}(M).\]
\end{theo*}
La preuve de ce r\'esultat utilise des arguments analogues \`a ceux
employ\'es par M.~Kashiwara \cite{kashiwara} dans le cas $p=1$.

Notons $b_L$ le polyn\^ome unitaire de plus bas degr\'e satisfaisant
l'identit\'e pr\'ec\'edente. La contribution de A.~Gyoja \cite{gyoja}
consiste dans le fait que $b_L$ est \`a racines dans $\Q_{<0}$.

\item[Etape 2]

Nous allons introduire deux nouvelles filtrations sur $M$.
\begin{enumerate}
\item
Soit $\sigma$ un c\^one convexe rationnel dans $\R_{\ge
  0}^p$. On note $\L(\sigma)$ l'ensemble des \'el\'ements primitifs du
$1$-squelette de $\sigma$ (i.e. $L \in \L(\sigma)$ si et seulement si
la droite engendr\'ee par $L$ est dans le $1$-squelette de $\sigma$ et
les coefficients de $L$ sont des entiers sans facteurs communs). Pour
tout $w \in \Z^p$, on note~:
\[{}^\sigma V_w(M)= \sum_{\{w'\in \Z^p \, | \forall L \in \L(\sigma)
  \, L(w')\le L(w)\}} V_{w'}(M).\]
(Voir la figure \ref{fig:Vsigma} o\`u $p=2$ et $\sigma$ est engendr\'e par
$L_1$ et $L_2$~: $m\in {}^\sigma V_w(M) \iff m$ est repr\'esent\'e par
un op\'erateur $P\in\Dnp$ dont le diagramme de Newton est dans le
quadrillage.)
\begin{figure}[h!]
\begin{center}
\includegraphics[angle=0, width=11cm]{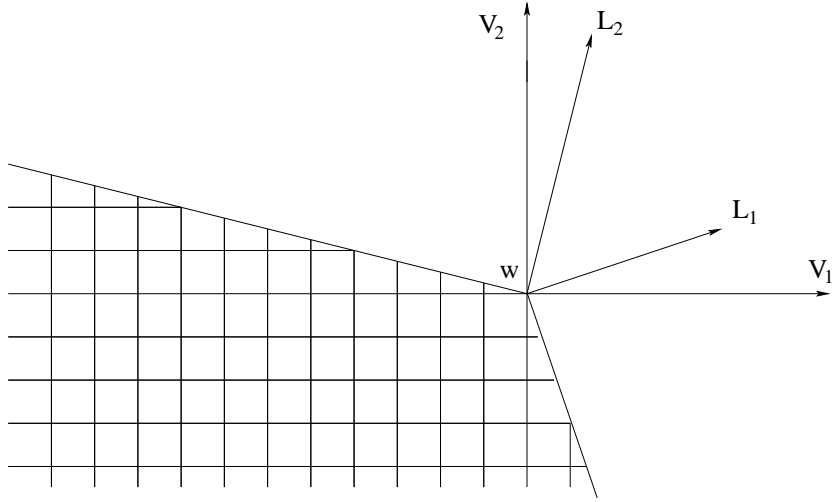}
\caption{${}^\sigma V_w(M)$}\label{fig:Vsigma}
\end{center}
\end{figure}

Pour tout $w\in \Z^p$, nous avons facilement les inclusions
suivantes~:
\[{}^\sigma V_w(M) \subseteq \bigcap_{L\in \sigma} V^L_{L(w)}(M)
\subseteq \bigcap_{L\in \L(\sigma)} V^L_{L(w)}(M),\]
la premi\`ere \'etant triviale et la seconde venant du fait que
$\L(\sigma)$ est inclus dans l'adh\'erence de~$\sigma$.

\begin{theo*}[S1](\cite{sabbah1} th. A.1.1 et prop. 2.2.1)
Il existe un \'eventail $\Sigma$ (dit \'eventail adapt\'e \`a $V(M)$)
constitu\'e de c\^ones poly\'edraux rationnels convexes tel que pour
tout c\^one $\sigma \in \Sigma$ et pour tout $w\in \Z^p$, on ait~:
\[{}^\sigma V_w(M) = \bigcap_{L\in \L(\sigma)} V^L_{L(w)}(M).\]
\end{theo*}
C'est ce th\'eor\`eme~S1 dont nous allons donner un \'enonc\'e et une
preuve plus \'el\'ementaires.
Disons un mot de la preuve. Dans l'appendice de \cite{sabbah1} en
collaboration avec F.J.~Castro Jim\'enez, C.~Sabbah d\'emontre
l'existence d'un \'eventail $\Sigma$ dit adapt\'e \`a la filtration
$V(M)$ apr\`es quoi il montre dans la proposition 2.2.1 que pour tout
c\^one d'un tel \'eventail, on a l'\'egalit\'e pr\'ec\'edente. La
preuve de l'existence d'un tel \'eventail m\'eriterait une mise au
point. En effet, elle s'appuie sur une division \`a param\`etre qui
conduirait \`a l'apparition de s\'eries formelles dans les variables
$\dxi$ de d\'erivation.
\item
Soit $\Vb(M)$ la filtration index\'ee par $\Z^p$ d\'efinie par~:
\[\Vb_w(M)=\bigcap_{L\in \UV} V^L_{L(w)}(M).\]
Comme cons\'equence du th\'eor\`eme~S1, nous avons~:

\begin{cor*}[S2]
\gauche
\begin{enumerate}
\item Pour tout $w \in \Z^p$,
\[\Vb_w(M)= \bigcap_{L\in \L(\Sigma)} V^L_{L(w)}(M),\]
o\`u $\L(\Sigma)$ d\'esigne l'ensemble des \'el\'ements primitifs du
$1$-squelette de $\Sigma$.
\item
Il existe $\kappa \in \N^p$ tel que pour tout $w\in \Z^p$, on ait
\[V_w(M) \subset \Vb_w(M) \subset V_{w+\kappa}(M).\]
\end{enumerate}
\end{cor*}

Faisons quelques commentaires sur la preuve de ce
r\'esultat. L'assertion (a) de S2 d\'ecoule trivialement de S1, en
effet~:
\begin{eqnarray*}
\Vb_w(M) & = & \bigcap_{\sigma \in \Sigma} \big(\bigcap_{L\in \sigma}
V^L_{L(w)}(M) \big)\\
 & = & \bigcap_{\sigma \in \Sigma} \big( \bigcap_{L\in \L(\sigma)}
V^L_{L(w)}(M)\big) \text{ par S1}\\
 & = & \bigcap_{L\in \L(\Sigma)} V^L_{L(w)}(M).
\end{eqnarray*}
En ce qui concerne, l'assertion (b), C.~Sabbah d\'emontre que si un
\'eventail $\Sigma$ satisfait l'\'enonc\'e du th\'eor\`eme S1 alors la
filtration $\Vb(M)$ est une bonne $V(\Dnp)$ filtration ce qui, par un
lemme usuel de comparaison entre bonnes filtrations fournit les
inclusions voulues (il faut noter que la premi\`ere inclusion de (b)
est trivialement vraie et que c'est la deuxi\`eme qui nous
int\'eresse).
\end{enumerate}
Finissons le rappel de la preuve. Notons $\delta$ la classe de $1$
dans le quotient $M=\Dnp/I$. Posons
\[b(s)= \prod_{L\in \L(\Sigma)} \Big(\prod_{-L(v+\kappa)<k\le 0}
b_L\big(L(s)-k\big)\Big)\]
Par les assertions (a) et (b) de S2, on constate que~:
\[b(-\ddt t) \delta \in \Vb_{-v-\kappa}(M) \subset V_{-v}(M),\]
ce qui signifie que $b(s) \in \B^v(f)$.

\end{description}

Nous allons maintenant \'enoncer les r\'esultats principaux du
pr\'esent article. 

On se donne un id\'eal $I$ de $\Dm$ (anneau des op\'erateurs
diff\'erentiels \`a coefficients analytiques faisant intervenir $m$
variables). Dans \cite{acg}, A.~Assi, F.~Castro~Jim\'enez et
M.~Granger ont introduit l'ensemble $\U$ des formes lin\'eaires $L$
pour lesquelles la filtration de $\Dm$ naturellement associ\'ee est
compatible avec la structure non commutative de $\Dm$ (cf. paragraphe
1.2). Ils ont \'etudi\'e le comportement du gradu\'e $\gr^L(I)$ lorque
$L$ varie dans $\U$. Consid\'erons la relation sur $\U$ telle que $L$
et $L'$ sont en relation si les gradu\'es $\gr^L(h(I))$ et
$\gr^{L'}(h(I))$ sont isomorphes ($h(I)$ d\'esigne l'homog\'en\'eis\'e
de $I$ dont nous rappelons la d\'efinition plus loin). Cette relation
est une relation d'\'equivalence donnant sur $\U$ une partition
constitu\'ee de c\^ones poly\'edraux rationnels convexes. L'eventail
ainsi constitu\'e est appel\'e \'eventail de Gr\"obner (analytique)
associ\'e \`a $h(I)$ et not\'e $\E(h(I))$.

Maintenenant, posons $m=n+p$ et reprenons les notations
pr\'ec\'edentes. Comme nous le verrons, on peut naturellement inclure
$\UV$ dans $\U$. On note alors $\EV=\EV(h(I))$ l'\'eventail obtenu
comme restriction de $\E(h(I))$ \`a $\UV$. Voici les r\'esultats qu'on
se propose de montrer ici.

\begin{theoi}\label{theo}
Pour tout c\^one $\sigma$ de $\EV$,
\[{}^\sigma V_w(M) = \bigcap_{L\in \L(\sigma)} V^L_{L(w)}(M).\]
\end{theoi}
L'analogue du corollaire S2 en d\'ecoule pour (a) et (b) en
rempla\c{c}ant le 1-squelette de $\Sigma$ par celui de $\EV$. En
effet, la preuve du corollaire S2 fonctionne pour n'importe quel
\'eventail d\`es que ce dernier satisfait le th\'eor\`eme S1.

\begin{theoi}\label{theo2}
Pour $p=2$, il existe $\kappa\in \N^2$ calculable \`a partir de bases
de Gr\"obner associ\'ees \`a chaque c\^one de l'\'eventail $\EV$ (ce
calcul sera d\'etaill\'e \`a la section 3) tel que pour tout $w\in
\Z^2$,
\[\Vb_{w+\kappa}(M) \subset V_w(M).\]
\end{theoi}
Pour $p\ge 3$, une g\'en\'eralisation de ce r\'esultat semble se
heurter \`a des difficult\'es techniques difficiles \`a r\'esoudre.

En r\'esum\'e, voici l'apport du pr\'esent article~:
\begin{itemize}
\item
Un \'enonc\'e et une preuve constructifs du th\'eor\`eme cl\'e S1 ce
qui fournit une approche plus \'el\'ementaire et plus constructive de
l'\'etape 2 de la preuve de C.~Sabbah et qui permet d'\'eviter la
notion d'\'eventail adapt\'e.
\item
Pour $p=2$, une preuve enti\`erement constructive de l'\'etape 2 dans
la preuve de C.~Sabbah.
\end{itemize}
Je signale que cet article constitue une partie de ma th\`ese
\cite{bahloul} dans laquelle on pourra trouver une autre preuve de
l'assertion (a) de S2 sans passer par le th\'eor\`eme \ref{theo}.\\
\\
D\'ecrivons le contenu de cet article.\\
Dans une premi\`ere section, nous faisons un certain nombre de
rappels concernant le th\'eor\`eme de division dans $\Dmz$ tel qu'il
est \'enonc\'e dans \cite{acg}, les notions de base standard et
d'\'eventail de Gr\"obner. La section 2 est consacr\'ee \`a la
d\'emonstration du th\'eor\`eme \ref{theo}. Dans la troisi\`eme
section, nous d\'emontrons le th\'eor\`eme \ref{theo2}.

\section{Rappels et r\'esultats pr\'eparatoires}

Dans les paragraphes qui suivent, nous rappelons, sans donner les
d\'emonstrations un certain nombre de notions et de r\'esultats qui
nous seront utiles dans les sections 2 et 3.

\subsection{Homog\'en\'eisation}

Comme dans \cite{castro-narvaez}, nous introduisons l'anneau
$\Dmz$. C'est avec cet anneau que les auteurs de \cite{acg} ont
introduit l'\'eventail de Gr\"obner analytique.

Dans ce paragraphe et jusqu'\`a 1.3.2, $x=(x_1,\ldots,x_m)$ et $\Dm$
d\'esigne l'anneau des op\'erateurs diff\'erentiels \`a coefficients
dans $\C\{x\}$. On d\'efinit l'anneau $\Dmz$ comme la $\C\{x\}$-alg\`ebre
engendr\'ee par $\dx{1}, \ldots, \dx{m}, z$ o\`u les seules relations
de commutation non triviales sont~:
\[[\dx{i}, c(x)]=\dxsur{c(x)}{i} z \text{ pour } i=1,\ldots,n \text{
  et } c(x) \in \C\{x\}.\]
Dans l'anneau $\Dm$, on note $\deg(P)$ le degr\'e total en les $\dxi$
d'un \'el\'ement $P$. Consid\'erons la filtration associ\'ee. Elle
s'\'etend naturellement \`a $\Dmz$ en consid\'erant le degr\'e total
en les $\dx{i}$ et $z$. Cette filtration fait de $\Dmz$ une alg\`ebre
gradu\'ee~:
\[\Dmz=\bigoplus_{d\in \N} \Dmz_d\quad  \text{ avec }
\quad \Dmz_d=\!\!\!\!\bigoplus_{k+|\beta|=d}
\C\{x\}\ddx^\beta z^k,\]
o\`u $\beta\in \N^m$, $\ddx^\beta=\dx{1}^{\beta_1} \cdots
\dx{m}^{\beta_m}$ et $|\beta|=\beta_1+\cdots+\beta_m$.\\
On dit qu'un op\'erateur $P\in \Dmz$ est homog\`ene (de degr\'e $d$)
si $P \in \Dmz_d$.\\
Pour $P\in \Dm$, on d\'efinit son homog\'en\'eis\'e $h(P)\in \Dmz$
comme suit.\\
On \'ecrit $P=\sum_{\beta}c_\beta(x) \ddx^\beta$ et on pose
$h(P)=\sum_\beta c_\beta(x)\ddx^\beta z^{d-|\beta|}$ o\`u $d=\deg(P)$,
ainsi $h(P)$ est homog\`ene de degr\'e $\deg(P)$.\\
De plus pour un id\'eal $I$ de $\Dm$, on d\'efinit $h(I)$ comme
l'id\'eal (\`a gauche) de $\Dmz$ engendr\'e par l'ensemble des
$h(P)$ avec $P$ parcourant $I$.

\subsection{Filtrations, divisions et bases standards}

Nous rappelons les notions usuelles de filtrations sur $\Dm$ et
$\Dmz$, le th\'eor\`eme de division dans $\Dmz$ de
\cite{acg}. Nous rappelons aussi les notions de bases standards et
bases standards r\'eduites minimales.

Soit $\U$ l'ensemble des formes lin\'eaires $L:\R^{2m} \to \R$,
$L(\alpha,\beta)=\sum_1^n e_i \alpha_i +\sum_1^n f_i \beta_i$ avec
pour tout $i=1,\ldots,n$, $e_i\le 0$ et $e_i+f_i\ge 0$. On \'etend
$\U$ \`a $\R^{2m+1}$ en posant $L(\alpha,\beta,k)=L(\alpha,\beta)$.\\
Pour $P\in \Dmz$ (resp. $P\in \Dm$), on \'ecrit $P=
\sum_{\alpha,\beta,k} a_{\alpha, \beta, k} x^\alpha \ddx^\beta z^k$
(avec $a_{\alpha, \beta, k}=0$ pour $k >0$ si $P\in \Dm$). On
d\'efinit le diagramme de Newton de $P$, not\'e $\DN(P) \subset
\N^{2m+1}$ (resp. $\subset \N^{2m}$), comme l'ensemble des
$(\alpha,\beta,k) \in \N^{2m+1}$ pour lesquels $a_{\alpha,\beta,k} \ne
0$.\\
Etant donn\'e $L \in \U$ et $P \in \Dmz$ (ou $P\in \Dm$), on d\'efinit
l'ordre de $P$ par rapport \`a $L$ comme $\ord^L(P)=\max
L(\DN(P))$. Cette ordre induit une filtration $F^L$ sur $\Dmz$ ou
$\Dm$ index\'ee par $L(\N^{2m+1})$ donn\'ee par~:
\[F^L_k(\Dmz)=\{P\in \Dmz; \, \ord^L(P) \le k\}\]
ainsi qu'un gradu\'e associ\'e $\dps\gr^L(\Dmz)= \bigoplus_{k\in
L(\N^{2m+1})} F^L_k(\Dmz)/F^L_{<k}(\Dmz)$.\\
Pour $P\in \Dmz$, on note
$\sigma^L(P)\in \gr^L(\Dnz)$ le symbole principal de $P$, i.e. la
classe de $P$ dans quotient $F^L_k(\Dmz)/F^L_{<k}(\Dmz)$ avec
$k=\ord^L(P)$. Si $J$ est un id\'eal de $\Dmz$, on a une filtration
$F^L(J)$ induite qui conduit \`a un id\'eal gradu\'e de $\gr^L(\Dmz)$
qu'on note $\gr^L(J)$ et qui est engendr\'e par l'ensemble des
$\sigma^L(P)$ pour $P\in J$.\\
\\
Pour une forme $L \in \U$, on d\'efinit deux ordres~: $<_L$ sur
$\N^{2m}$ et $<_L^h$ sur $\N^{2m+1}$~:
\[(\alpha,\beta) <_L (\alpha',\beta') \iff 
\begin{cases}
L(\alpha,\beta) < L(\alpha',\beta') \\
\text{ou } \big(  \text{ \'egalit\'e et } |\beta|<|\beta'| \big) \\
\text{ou } \big( \text{ \'egalit\'es et } (\alpha,\beta) >_0
(\alpha',\beta') \big),
\end{cases}\]
o\`u $<_0$ est un bon ordre compatible avec l'addition fix\'e pour
toute la suite.

\[(\alpha,\beta,k) <_L^h (\alpha',\beta',k') \iff
\begin{cases}
k+|\beta|<k'+|\beta'| \\
\text{ou } \big( \text{ \'egalit\'e et } (\alpha,\beta) <_L
(\alpha',\beta') \big).
\end{cases}\]

Pour $P\in \Dmz$, on d\'efinit l'exposant privil\'egi\'e
$\exp_{<_L^h}(P)= \max _{<_L^h}(\DN(P))$ et le mon\^ome privil\'egi\'e
$\mp_{<_L^h}(P)=(x,\ddx,z)^{\exp_{<_L^h}(P)}$ . De m\^eme pour $P \in
\Dm$ et $<_L$ \`a la place de $<_L^h$. Notons que
$\exp_\prec(PQ)=\exp_\prec(P) +\exp_\prec(Q)$ d\`es que $\prec$ est
compatible avec l'addition ce qui est le cas pour $<_L$ et $<_L^h$.
Rappelons le th\'eor\`eme de division dans $\Dmz$ de \cite{acg}.

Soit $L \in \U$. Soient $Q_1,\ldots,Q_r$ une famille
d'\'el\'ements de $\Dmz$. Notons $(\Delta_1,\ldots,
\Delta_r,\bar{\Delta})$ la partition de $\N^{2m+1}$ d\'efinie \`a
partir des $\exp_{<_L^h}(Q_j)$ :
\begin{itemize}
\item $\Delta_1=\exp_{<_L^h}(Q_1)+\N^{2m+1}$
\item $\Delta_j=(\exp_{<_L^h}(Q_j)+\N^{2m+1}) \smallsetminus
  (\bigcup_{i=1}^{i=j-1} \Delta_i)$ pour $j=2,\ldots,r$
\item $\bar{\Delta}=\N^{2m+1} \smallsetminus (\bigcup_{j=1}^{j=r}
  \Delta_j)$
\end{itemize}
\begin{theo} (\cite{acg} Theorem 7)\label{theo:div}
Pour tout $P\in \Dmz$ il existe un unique $(q_1,\ldots,q_r,R) \in
(\Dmz)^{r+1}$ tel que :
\begin{enumerate}
\item $P=q_1 Q_1 + \cdots + q_r Q_r+R$
\item pour tout $j=1,\ldots,r$, si $q_j \ne 0$ alors
  $\DN(q_j)+\exp_{<_L^h}(Q_j) \subset \Delta_j$
\item si $R \ne 0$ alors $\DN(R) \subset \bar{\Delta}$
\end{enumerate}
On appelle $R$ le reste de la division de $P$ par les $Q_j$
relativement \`a $<_L^h$.
\end{theo}

\begin{cor}\gauche
\begin{itemize}
\item
$\exp_{<_L^h}(P)=\max_{<_L^h}\{\exp_{<_L^h}(q_jQ_j),\, j=1,\ldots,r;
\, \exp_{<_L^h}(R)\}$
\item
$\ord^L(P)=\max\{\ord^L(q_jQ_j),\, j=1,\ldots,r;\, \ord^L(R)\}$
\end{itemize}
\end{cor}
Soit $J$ un id\'eal de $\Dmz$ et $Q_1,\ldots,Q_r \in J$. On dit que
$Q_1,\ldots,Q_r$ forment une $<_L^h$-base standard de $J$ si pour tout
$P$ dans $J$, le reste de la division de $P$ par $Q_1,\ldots,Q_r$ est
nul. Consid\'erons l'ensemble des exposants de $J$,
$\Exp_{<_L^h}(J)=\{\exp_{<_L^h}(P), \, P\in J\smallsetminus
0\}$. Etant donn\'es $Q_1,\ldots,Q_r \in J$, les deux assertions
suivantes sont \'equivalentes (gr\^ace aux th\'eor\`eme de division)~:
\begin{itemize}
\item Les $Q_j$ forment une $<_L^h$-base standard de $J$
\item $\Exp_{<_L^h}(J)=\bigcup_{j=1}^r (\exp_{<_L^h}(Q_j)+\N^{2m+1})$
\end{itemize}
L'existence d'une base standard est assur\'ee par le lemme de Dickson
qui dit que si une partie $E$ de $\N^q$ ($q\in \N\smallsetminus 0$)
v\'erifie $E+\N^q=E$ (ce qui est le cas de $\Exp_{<_L^h}(J)$) alors il
existe $F\subset E$ fini tel que $E=\cup_{e\in F} (e+\N^q)$.

Il est facile de voir qu'une base standard n'est pas unique, c'est
pour cela qu'il existe une notion de base standard r\'eduite
minimale~:
\begin{defin*}
Soit $Q_1,\ldots,Q_r$ une $<_L^h$-base standard de $J\subset \Dmz$ et
soit $e_j=\exp_{<_L^h}(Q_j)$ pour $j=1,\ldots,r$.
\begin{itemize}
\item
On dit qu'elle est minimale si pour toute partie finie $F$ de
$\N^{2m+1}$ on a l'implication suivante~:
\[\Exp_{<_L^h}(J)=\bigcup_{e\in F}(e+\N^{2m+1}) \Rightarrow
\{e_1,\ldots,e_r\} \subseteq F.\]
\item
On dit qu'elle est r\'eduite si les $Q_j$ sont unitaires (i.e. le
coefficient du mon\^one privil\'egi\'e est $1$) et si pour tout $j$~:
\[(\DN(Q_j) \smallsetminus e_j) \subset (\N^{2m+1} \smallsetminus
\Exp_{<_L^h}(J)).\]
\end{itemize}
Il existe une unique $<_L^h$-base standard r\'eduite minimale d'un
id\'eal $J\subset \Dmz$.
\end{defin*}
Notons que si l'id\'eal $J$ est homog\`ene alors la base standard
r\'eduite minimale est constitu\'ee d'\'el\'ements homog\`enes.
Nous terminons ce paragraphe par un lemme qui illustre l'int\'er\^et
des bases standards r\'eduites minimales. Nous l'appliquerons dans la
section suivante.

\begin{lem}\label{lem:utile}
Soient $<_1$ et $<_2$ deux ordres sur $\N^{2m+1}$ qui permettent de
faire des divisions dans $\Dmz$ (par exemple $<_i=<_{L_i}^h$ avec
$L_1$ et $L_2$ deux formes). Soit $\{Q_1,\ldots,Q_r\}$ la base
standard r\'eduite minimale d'un id\'eal $J \subset \Dmz$ par rapport
\`a l'ordre $<_1$. Supposons que pour tout $j$, $\exp_{<_1}(Q_j)$ et
$\exp_{<_2}(Q_j)$ soient \'egaux. Alors $\{Q_1,\ldots,Q_r\}$ est aussi
la base standard r\'eduite minimale de $J$ par rapport \`a $<_2$.
\end{lem}
Nous omettons la preuve. Disons simplement que ce qui importe lors
d'une division, c'est l'ensemble des $\exp_<(Q_j)$ et pas l'ordre
$<$. Ainsi, diviser par les $Q_j$ par rapport \`a $<_1$ ou $<_2$
donnera les m\^emes quotients et le m\^eme reste.


\subsection{Eventail de Gr\"obner}

Dans ce paragraphe, nous rappelons le r\'esultat principal de
\cite{acg} d\'ecrivant l'\'eventail de Gr\"obner d'un id\'eal. Nous
introduisons aussi le $V$-\'eventail de Gr\"obner $\EV(h(I))$, ce sera
l'objet principal des sections $2$ et $3$.

\subsubsection{}

On se donne un id\'eal $I$ de $\Dm$ et on consid\`ere $h(I) \subset
\Dmz$ son homog\'en\'eis\'e.
Pour $L$ et $L'$ dans $\U$, on d\'efinit la relation~:
\[L \sim L' \iff \gr^L(\Dmz) \simeq \gr^{L'}(\Dmz) \text{ et }
\gr^L(h(I)) \simeq \gr^{L'}(h(I)).\]
Cette relation est une relation d'\'equivalence sur $\U$.

\begin{theo}\label{theo:event_Grob}\cite{acg}
La partition sur $\U$ donn\'ee par la relation d'\'equivalence
pr\'ec\'edente est constitu\'ee de c\^ones poly\'edraux rationnels
convexes. Cette partition qu'on note $\E(h(I))$ est appel\'ee
\'eventail de Gr\"obner (analytique) de $h(I)$.\\
De plus pour chaque c\^one $\sigma \in \E(h(I))$, il existe
$Q_1,\ldots,Q_r \in h(I)$ homog\`enes tels que :
\begin{itemize}
\item pour tout $L,L' \in \sigma$, $\sigma^L(Q_j)=\sigma^{L'}(Q_j)$
  et $\exp_{<_L^h}(Q_j)=\exp_{<_{L'}^h}(Q_j)$ pour tout $j$.
\item pour tout $L\in \sigma$, l'ensemble $\{Q_1,\ldots,Q_r\}$ est la
  base standard r\'eduite minimale de $h(I)$ par rapport \`a $<_L^h$.
\end{itemize}
\end{theo}
Remarquons que la deuxi\`eme condition nous montre que sur un c\^one
$\sigma$, l'ensemble $\Exp_{<_L^h}(h(I))$ des exposants de $h(I)$ par
rapport \`a $<_L^h$ est constant lorsque $L$ parcourt $\sigma$.\\

\subsubsection{}

A partir de maintenant, on travaillera dans $\Dnp$ et $\Dnpz$
(i.e. $m=n+p$). Pour $j=1,\ldots,p$, on notera $V_j \in \U$ la forme
lin\'eaire donn\'ee par~: $V_j(\alpha,\mu,\beta,\nu)=\nu_j -\mu_j$
o\`u $\alpha,\beta \in \N^n$ et $\mu,\nu\in \N^p$. Cette forme
lin\'eaire donne lieu \`a une filtration qu'on note aussi $V_j$ et qui
n'est rien d'autre que la $V$-filtration de Kashiwara-Malgrange
associ\'ee \`a la variable $t_j$ (rappelons que dans $\Dnp$, les
variables sont $x=(x_1,\ldots,x_n)$ et $t=(t_1,\ldots, t_p)$). On note
alors $V$ la multifiltration $V=(V_1,\ldots,V_p)$.\\
Nous noterons $\UV\subset \U$ le sous-ensemble des formes lin\'eaires
$L$ s'\'ecrivant~:
\[L=l_1 V_1+\cdots+l_p V_p,\]
avec $(l_1,\ldots,l_p)\in (\R_{\ge 0})^p$, aussi nous identifierons
$\UV$ et $(\R_{\ge 0})^p$.

A partir de maintenant, pour $L\in \UV$, on notera $V^L$ la filtration
associ\'ee \`a $L$ (conform\'ement aux notations de \cite{sabbah1}).
Notons que pour $L \in \UV$, $\gr^L(\Dnpz)$ est isomorphe \`a un
sous-anneau de $\Dnpz$, ainsi nous consid\'ererons que tous les
calculs se feront dans $\Dnpz$. Maintenant, consid\'erons la
restriction de $\E(h(I))$ \`a $\UV$. Ceci donne lieu au $V$-\'eventail
de Gr\"obner de $h(I)$ que l'on note $\EV(h(I))$. Comme cons\'equence,
on obtient un analogue du th\'eor\`eme \ref{theo:event_Grob}~:
\begin{cor}\label{cor:V_event_Grob}
Pour chaque c\^one $\sigma$ de $\EV(h(I))$, il existe
$Q_1,\ldots,Q_r \in h(I)$ homog\`enes tels que~:
\begin{itemize}
\item pour toute $L,L' \in \sigma$, $\sigma^L(Q_j)=\sigma^{L'}(Q_j)$
  et $\exp_{<_L^h}(Q_j)=\exp_{<_{L'}^h}(Q_j)$ pour tout $j$.
\item pour tout $L\in \sigma$, l'ensemble $\{Q_1,\ldots,Q_r\}$ est la
  base standard r\'eduite minimale de $h(I)$ par rapport \`a $<_L^h$.
\end{itemize}
\end{cor}
On voit que pour chaque c\^one $\sigma$, il existe des
$Q_j$ qui, pour tout $L$ de $\sigma$, forment la base standard
r\'eduite minimale de $h(I)$ par rapport \`a $<_L^h$. S'il n'y pas de
confusion possible, on l'appelle {\bf la base standard de $h(I)$
associ\'ee \`a $\sigma$}.

Terminons ce paragraphe par quelques remarques et notations.
Un c\^one $\sigma$ de
$\EV(h(I))$ n'est pas n\'ecessairement ouvert (il peut \^etre
"semi-ouvert"). Voici comment on d\'efinit $\L(\sigma)$~:\\
on consid\`ere d'abord l'adh\'erence $\bar{\sigma}$ de $\sigma$.
Il existe $L_1,\ldots,L_q \in \UV$ qu'on suppose primitifs (avec $q
\ge 1$ pouvant \^etre plus grand que $p$) tels que~:
\[\bar{\sigma}=\{L=r_1 L_1+\cdots r_q L_q; \, r_i \ge 0\}.\]
Supposons $q$ minimal alors $\L(\sigma)$ est l'ensemble
$\{L_1,\ldots,L_q\}$. L'ensemble $\L(\EV(h(I)))$ n'est alors rien
d'autre que la r\'eunion des $\L(\sigma)$ avec $\sigma \in \EV(h(I))$
(remarquons qu'il est constitu\'e d'\'el\'ements entiers car tout
c\^one $\sigma$ est rationnel).

On d\'efinit l'int\'erieur de $\sigma$ comme l'ensemble des
combinaisons \`a coefficients strictement positifs des $L_i$. On note
$\rangle L_1,\ldots, L_q \langle$ le c\^one ouvert engendr\'e par les
$L_i$. D'autre part si $L_1$ et $L_2$ sont dans $\UV$, on note
$\rangle L_1,L_2 \rangle=\{r_1L_1 +r_2L_2; \, r_1>0, \, r_2\ge 0\}$ le
c\^one semi-ouvert contenant $L_2$ mais non $L_1$.

Dans la suite, on \'ecrira $\EV$ \`a la place de $\EV(h(I))$.

\section{D\'emonstration du th\'eor\`eme \ref{theo}}

Soit $\sigma$ un c\^one de $\EV$. Nous avons vu dans le corollaire
\ref{cor:V_event_Grob} qu'il existe une famille $Q_1,\ldots,Q_r$ de
$h(I)$ qui est la base standard r\'eduite minimale de $h(I)$ par
rapport \`a $<_L^h$ et ce pour tout $L$ dans $\sigma$. Maintenant, que
se passe-t-il pour une forme lin\'eaire $L$ appartenant \`a
l'adh\'erence de $\sigma$ mais qui n'est pas dans $\sigma$ (ce qui
peut arriver pour un \'el\'ement de $\mathcal{L}(\sigma)$) ? Bien
entendu les $Q_j$ ne forment pas n\'ecessairement une base standard de
$h(I)$ pour l'ordre $<_L^h$. Par contre il est possible, et c'est le
but du r\'esultat suivant, de construire un ordre que nous noterons
$\lhd_L^\sigma$ (car il d\'epend de $L$ et de $\sigma$) pour lequel
les $Q_j$ forment la base standard r\'eduite minimale de $h(I)$. Voici
comment on d\'efinit l'ordre en question.\\
D'abord on fixe une forme lin\'eaire $L_\sigma$ dans l'int\'erieur de
$\sigma$ et pour $(\alpha,\mu,\beta,\nu,k)$ et
$(\alpha',\mu',\beta',\nu',k')$ dans $\N^{n+p+n+p+1}$ on pose~:
\[(\alpha,\mu,\beta,\nu,k) \lhd_L^\sigma (\alpha',\mu',\beta',\nu',k)
\iff\]
\[\begin{cases}
k+|\beta|+|\nu|<k'+|\beta'|+|\nu'| \\
\text{ou } \big( = \text{ et } L(\alpha,\mu,\beta,\nu)<
L(\alpha',\mu',\beta',\nu') \big) \\
\text{ou } \big( = \text{ et } = \text{ et } (\alpha,\mu,\beta,\nu)
<_{L_\sigma} (\alpha',\mu',\beta',\nu') \big).
\end{cases}\]

\begin{prop}
Soit $L \in \UV$ appartenant \`a l'adh\'erence de $\sigma \in \EV$ et
soit $Q_1,\ldots,Q_r$ la base standard de $h(I)$ associ\'ee au
c\^one $\sigma$. Alors pour tout $j=1,\ldots,r$ et pour tout $L'$ dans
$\rangle L, L_\sigma \rangle$~:
\[\exp_{<_{L'}^h}(Q_j)=\exp_{\lhd_L^\sigma}(Q_j).\]
\end{prop}
Comme {\bf cons\'equence} de cette proposition, nous obtenons que {\bf
  les $Q_j$ forment la base standard r\'eduite minimale de $h(I)$ pour
  l'ordre $\lhd_L^\sigma$}. En effet, il suffit pour cela d'appliquer
le lemme \ref{lem:utile}.

\begin{figure}[h!]
\begin{center}
\includegraphics[angle=0, width=11cm]{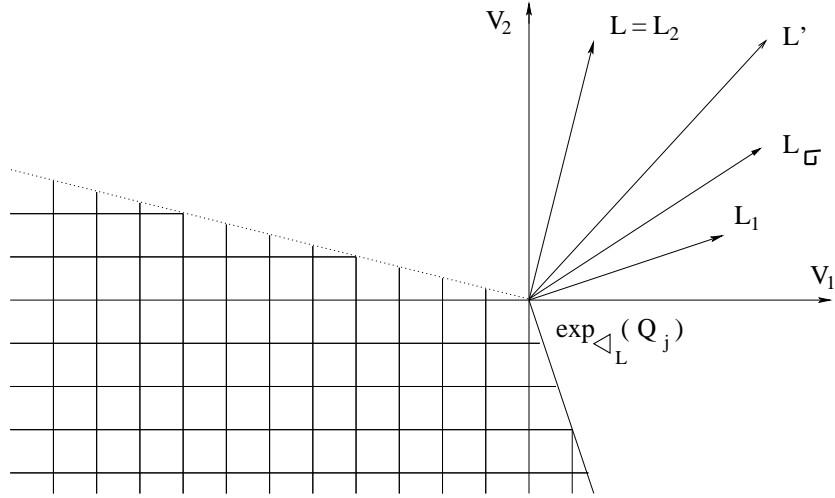}
\caption{Diagramme de Newton d'un $Q_j$ associ\'e \`a $\sigma=\rangle
  L_2, L_1 \rangle$}\label{fig:Qj_sans_pentes}
\end{center}
\end{figure}

\begin{proof}[D\'emonstration]
En utilisant le corollaire \ref{cor:V_event_Grob}, on obtient
\begin{eqnarray*}
\exp_{\lhd_L^\sigma}(Q_j) & = & \exp_{\lhd_L^\sigma}(\sigma^L(Q_j)) \\
 & = & \exp_{<_{L_\sigma}^h}(\sigma^L(Q_j)) \text{ par d\'efinition de
   } \lhd_L^\sigma \text{ et } <_{L_\sigma}^h\\
 & = & \exp_{<_{L_\sigma}^h}(\sigma^{L_\sigma}( \sigma^L(Q_j))) \\
 & = & \exp_{<_{L_\sigma}^h}(\sigma^{L_\sigma}(Q_j)) \text{ car }
 \sigma^{L''}(Q_j) \text{ est constant } \forall L''\in \rangle L,
 L_\sigma \rangle,\\
 & & \text{voir la figure \ref{fig:Qj_sans_pentes} ci-apr\`es.}\\
 & = & \exp_{<_{L_\sigma}^h}(Q_j)\\
 & = & \exp_{<_{L'}^h}(Q_j).
\end{eqnarray*}
\end{proof}

Rappelons que $I$ est un id\'eal de $\Dnp$, $M=\Dnp/I$ et $\delta$ est
la classe de $1$ dans $M$.
Soit $\sigma$ un c\^one de $\EV$. Notons $L_1,\ldots,L_q$ les
\'el\'ements de $\mathcal{L}(\sigma)$.
\begin{lem}\label{lem:cle}
Soit $i_0 \in \{1,\ldots,q\}$, $m \in V^{L_{i_0}}_{\lambda_{i_0}}(M)$
avec $\lambda_{i_0} \in \Q$ et $P \in \Dnp$ tel que $P \delta=m$ et
$\ord^{L_{i_0}}(P)>\lambda_{i_0}$ alors il existe $P' \in \Dnp$ tel
que :
\begin{itemize}
\item $P-P' \in I$ c'est-\`a-dire $P' \delta=m$
\item $\ord^{L_{i_0}}(P')<\ord^{L_{i_0}}(P)$
\item $\ord^{L_i}(P') \le \ord^{L_i}(P)$ pour $i\in \{1,\ldots,q\}
  \smallsetminus \{i_0\}$.
\end{itemize}
\end{lem}
En d'autres termes, il est possible de faire baisser l'ordre par
rapport \`a l'un des $L_i$ sans augmenter l'ordre par rapport aux
autres $L_i$. Avec ce lemme, nous sommes en mesure de donner une

\begin{proof}[D\'emonstration du th\'eor\`eme \ref{theo}]
Soit $m \in \bigcap_{L \in \mathcal{L}(\sigma)}
V^L_{L(v)}(M)$ alors pour $i=1,\ldots,q$, il existe $P_i \in
\Dnp$ tel que $P_i \delta=m$ et $\ord^{L_i}(P_i) \le L_i(v)$. 
On pose $\tilde{P}_1=P_1$. En appliquant un nombre fini de fois
le lemme avec $i_0=2$ (la premi\`ere fois avec $P=\tilde{P}_1$ et
$\lambda_{i_0}=\ord^{L_{i_0}}(\tilde{P}_1)$),
on construit $\tilde{P}_2$ tel que $\tilde{P}_2\delta=m$ et
$\ord^{L_i}(\tilde{P}_2) \le L_i(v)$ pour $i=1,2$.
On recommence le processus avec $i_0=3$ et $P=\tilde{P}_2$,
etc. Apr\`es un nombre fini d'\'etapes on obtient $\tilde{P}_q \in \Dnp$
tel que $\tilde{P}_q \delta=m$ et pour tout $i=1,\ldots,q$,
$\ord^{L_i}(\tilde{P}_q) \le L_i(v)$. Ceci montre bien que $m \in
{}^{\sigma}V_v(M)$.
\end{proof}

Pour finir cette section, il ne reste plus qu'\`a d\'emontrer le lemme
pr\'ec\'edent~:
\begin{proof}[D\'emonstration de \ref{lem:cle}]
Pour simplifier les notations, nous ferons la d\'emonstration avec
$i_0=1$, $\lambda=\lambda_{i_0}$ et $\lhd_{L_i}= \lhd_{L_i}^\sigma$.
Notons $Q_1,\ldots, Q_r$ la base standard de $h(I)$ associ\'ee au
c\^one $\sigma$.

Par hypoth\`ese, il existe $P_1 \in \Dnp$ avec $P_1 \delta=m$ et
$\ord^{L_1}(P_1)\le \lambda$. Il existe $l_0,l,l_1 \in \N$ tels que
$z^{l_0}h(P-P_1)=z^lh(P)-z^{l_1}h(P_1)$. On pose alors $H=z^lh(P)$,
$H_1=z^{l_1}h(P_1)$ et $H_0=H-H_1$ et comme $P-P_1 \in I$, on a
$H_0 \in h(I)$.\\

Consid\'erons la division de $H_0$ par les $Q_j$ relativement \`a
l'ordre $\lhd_{L_1}$~:
\[H_0=\sum_{j=1}^r q_j Q_j \text{ avec }
\mathcal{N}(q_j)+\exp_{\lhd_{L_1}}(Q_j) \subset \Delta_j \text{ pour
  tout } j\]
o\`u les $\Delta_j \subset \N^{2n+2p+1}$ forment la partition de
$\Exp_{\lhd_{L_1}}(h(I))$ associ\'ee aux exposants privil\'egi\'es des
$Q_j$ (voir le th\'eor\`eme de division \ref{theo:div} rappel\'e \`a
  la section $1$).\\
Comme pour tout $i,j$, les exposants $\exp_{\lhd_{L_1}}(Q_j)$ et
$\exp_{\lhd_{L_i}}(Q_j)$ sont \'egaux, la division pr\'ec\'edente est
aussi une division relativement aux ordres $\lhd_{L_2}, \ldots,
\lhd_{L_q}$. Par cons\'equent, pour tout $i=1,\ldots,q$ et
$j=1,\ldots,r$~:
\[\ord^{L_i}(H_0) \ge \ord^{L_i}(q_j Q_j).\]
Notons $J$ l'ensemble des $j\in \{1,\ldots,r\}$ pour lesquels
$\ord^{L_1}(H_0)=\ord^{L_1}(q_j Q_j)$, nous avons alors~:
\[\sigma^{L_1}(H_0)=\sigma^{L_1}(H)= \sum_{j\in J} \sigma^{L_1}(q_j)
\sigma^{L_1}(Q_j).\]
On consid\`ere et on note $\dps W=\sum_{j\in J}
\sigma^{L_1}(q_j)Q_j$. C'est un \'el\'ement de $h(I)$.\\
Posons $H'=H-W$. Nous allons montrer les deux assertions suivantes.
\begin{enumerate}
\item $\ord^{L_1}(H') < \ord^{L_1}(H),$
\item $\ord^{L_i}(H') \le \ord^{L_i}(H)$ pour $i=2,\ldots,q$.
\end{enumerate}

\begin{enumerate}
\item
On a clairement $\sigma^{L_1}(H)=\sigma^{L_1}(W)$. Par cons\'equent,
\[H'=(H-\sigma^{L_1}(H)) - (W -\sigma^{L_1}(W)).\]
On voit alors facilement que les deux termes entre parenth\`eses ont
un $L_1$-ordre strictement inf\'erieur \`a celui de $H$.

\item
Fixons $i$ entre $2$ et $q$.\\
En utilisant \ref{cor:V_event_Grob}, pour tout $j$, on a~:
\begin{equation}\label{ident1}
\exp_{\lhd_{L_i}}(\sigma^{L_1}(Q_j))= \exp_{\lhd_{L_i}}(Q_j).
\end{equation}
D'autre part, par construction des $q_j$, pour tout $j$ on a~:
$\DN(q_j)+ \exp_{\lhd_{L_i}}(Q_j)\subset \Delta_j$, ce qui entraine
ceci pour tout $j$~:
\begin{equation}\label{ident2}
\DN(\sigma^{L_1}(q_j))+\exp_{\lhd_{L_i}}(\sigma^{L_1}(Q_j))\subset
\Delta_j.
\end{equation}
Maintenant, nous avons $\sigma^{L_1}(W)=\sum_{j\in J} \sigma^{L_1}(q_j)
\sigma^{L_1}(Q_j)$. Par la relation (\ref{ident2}), on peut dire que
cette \'ecriture est le r\'esultat de la division de
$\sigma^{L_1}(W)$ par $\{\sigma^{L_1}(Q_j), \, j\in J\}$ relativement
\`a $\lhd_{L_i}$, par cons\'equent~:
\begin{equation}\label{ident3}
\exp_{\lhd_{L_i}}(\sigma^{L_1}(W)) =\max_{j\in J}
\{\exp_{\lhd_{L_i}} (\sigma^{L_1}(q_j) \sigma^{L_1}(Q_j)) \}.
\end{equation}
De la m\^eme mani\`ere, on peut montrer que
\begin{equation}\label{ident4}
\exp_{\lhd_{L_i}}(W) =\max_{j\in J}
\{\exp_{\lhd_{L_i}} (\sigma^{L_1}(q_j) Q_j) \}.
\end{equation}
Par cons\'equent, gr\^ace \`a (\ref{ident3}), (\ref{ident4}) et
(\ref{ident1}), on obtient l'\'egalit\'e
$\exp_{\lhd_{L_i}}(\sigma^{L_1}(W)) =\exp_{\lhd_{L_i}}(W)$, ce qui
donne en particulier $\ord^{L_i}(W)=\ord^{L_i}(\sigma^{L_1}(W))$.
D'o\`u les \'egalit\'es et in\'egalit\'es suivantes~:
\begin{eqnarray*}
\ord^{L_i}(W) & = & \ord^{L_i}(\sigma^{L_1}(W)) \\
 & = & \ord^{L_i}(\sigma^{L_1}(H)) \text{ car }
 \sigma^{L_1}(W)=\sigma^{L_1}(H) \\
 & \le & \ord^{L_i}(H). 
\end{eqnarray*}
Par cons\'equent, $\ord^{L_i}(H') \le \ord^{L_i}(H)$.
Les deux assertions sont d\'emontr\'ees.
\end{enumerate}
Maintenant, effectuons la sp\'ecialisation $z=1$ (qui est un morphisme
d'alg\`ebres $\Dnpz \to \Dnp$) et posons $P'=H'_{|z=1}=P-W_{|z=1}$.
On a $W\in h(I)$ donc $W_{|z=1} \in I$ et $P'\delta=m$. Apr\`es
sp\'ecialisation, les assertions $1$ et $2$ deviennent~:
$\ord^{L_i}(P') \le \ord^{L_i}(P)$ pour tout $i=1,\ldots,q$ avec
in\'egalit\'e stricte pour $i=1$. Le lemme est d\'emontr\'e.
\end{proof}

\section{D\'emonstration du th\'eor\`eme \ref{theo2}}

Dans cette section, nous allons d\'emontrer le th\'eor\`eme
\ref{theo2}. Nous allons dans un premier temps \'enoncer plus
pr\'ecis\'ement l'assertion en question, puis nous en \'ecrirons la
preuve. Nous verrons qu'elle consiste essentiellement en une analyse
raffin\'ee du lemme pr\'ec\'edent \ref{lem:cle}, en effet nous ferons
ce qu'on pourrait appeler un ``contr\^ole de la mont\'ee de
l'ordre'' par rapport \`a la forme $V_1$.\\
Rappelons que dans ce paragraphe, $p$ \'egale $2$.

\begin{nota*}
\gauche
\begin{itemize}
\item
Soient $L_1,L_2$ deux formes non nulles de $\UV$. Ecrivons
$L_i=a_iV_1+b_iV_2$ avec $a_i,b_i \ge 0$. On dit que $L_1$ est
inf\'erieure (resp. strictement inf\'erieure) \`a $L_2$ si $b_1/a_1
\le b_2/a_2$ (resp. $b_1/a_1 < b_2/a_2$). On abr\`egera cette notion
en notant $L_1 \le L_2$ (resp. $L_1 < L_2$). Par convention
$b/0=+\infty$, toute forme $L$ est inf\'erieure \`a $V_2$.
\item
Soient $L_1 \ne L_2$ dans $\UV$ et $H\in \Dnpz$. On dit que $H$ est
$L_1$-homog\`ene si $H=\sigma^{L_1}(H)$. On dit que $H$ est
$(L_1,L_2)$-homog\`ene si $H=\sigma^{L_1}(\sigma^{L_2}(H))$.
\item
Soit $L$ une forme dans $\UV$. Nous noterons $\lhd_L$ l'ordre sur
$\N^{n+2+n+2+1}$ donn\'e par :
\[(\alpha,\mu,\beta,\nu,k) \lhd_L (\alpha',\mu',\beta',\nu',k) \iff\]
\[\begin{cases}
k+|\beta+\nu|<k'+|\beta'+\nu'| \\
\text{ou } \big( = \text{ et } L(\alpha,\mu,\beta,\nu)<
L(\alpha',\mu',\beta',\nu') \big) \\
\text{ou } \big( = \text{ et } = \text{ et } (\alpha,\mu,\beta,\nu)
<_{V_1} (\alpha',\mu',\beta',\nu') \big).
\end{cases}\]
Nous remarquons qu'en adoptant les notations du paragraphe
pr\'ec\'edent et en posant $\sigma=\rangle V_1, L \langle$ (avec $L\ne
V_1$) alors on a : $\lhd_L=\lhd_L^\sigma$. Si par contre $L=V_1$ alors
$\lhd_L=<_{V_1}^h$.
\end{itemize}
\end{nota*}

\begin{figure}[h!]
\begin{center}
\includegraphics[angle=0, width=9cm]{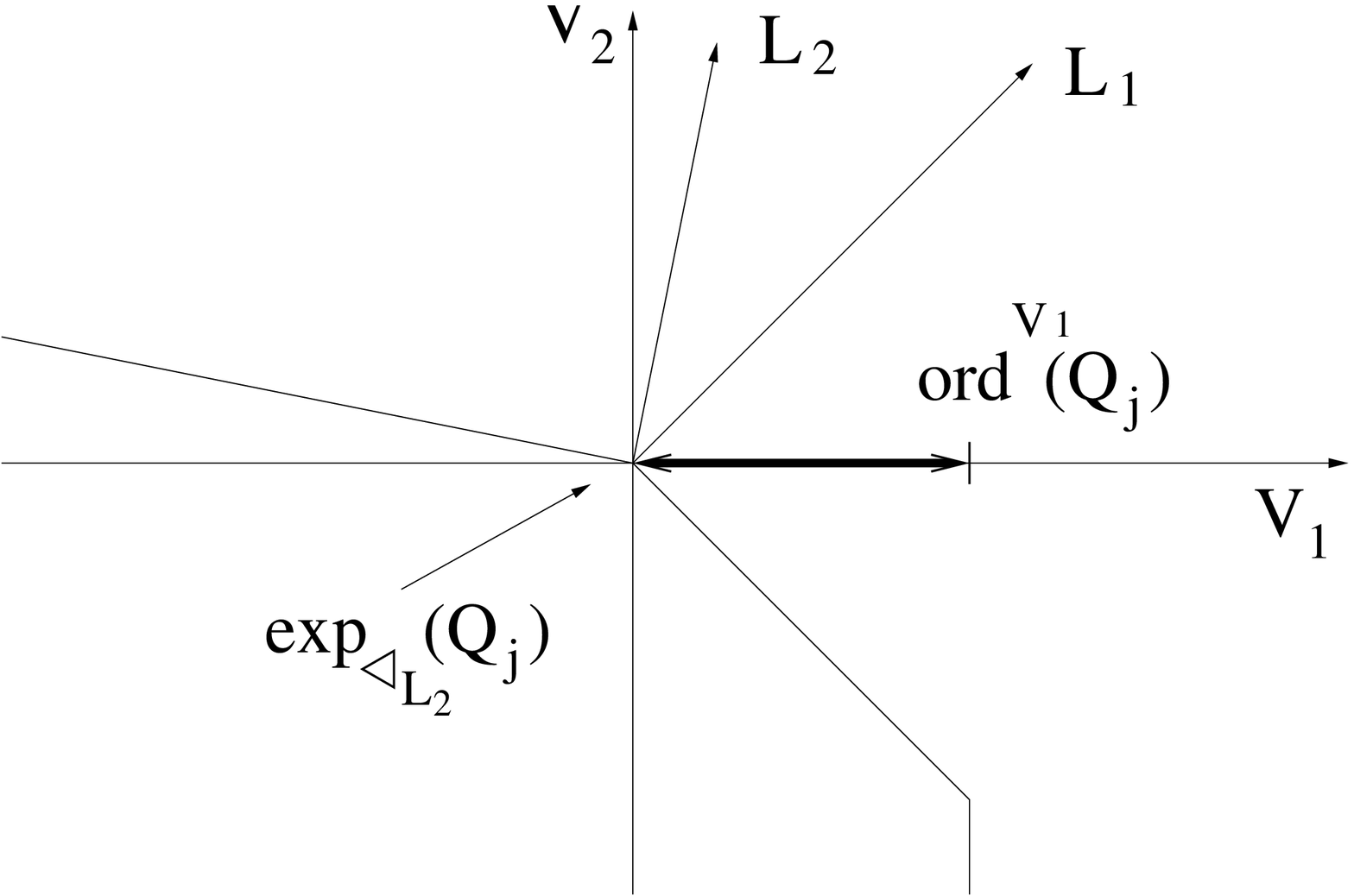}
\caption{$\ord^{V_1}(Q_j)-
  \ord^{V_1}(\exp_{\lhd_{L_2}}(Q_j))$}\label{fig:kappa1}
\end{center}
\end{figure}

Soit $\sigma$ un c\^one de $\EV$ de dimension $2$ (maximale) et
$\{L_1,L_2\}=\mathcal{L}(\sigma)$ avec $L_1 <L_2$. Notons
$Q_1,\ldots,Q_r$ la base standard de $h(I)$ associ\'ee \`a
$\sigma$. On d\'efinit $\kappa_\sigma^1 \in \N$ par~:
\[\kappa_\sigma^1=\max\{ \ord^{V_1}(Q_j)-\ord^{V_1}(\sigma^{L_2}(Q_j))
, \, j=1,\dots,r\}.\]
Avec les notations pr\'ec\'edentes, nous avons (voir la figure
\ref{fig:kappa1})~:
\[\ord^{V_1}(\sigma^{L_2}(Q_j))=\ord^{V_1}(\exp_{\lhd_{L_2}}(Q_j)).\]

Maintenant on d\'efinit $\kappa^1 \in \N$ comme le maximum des
$\kappa_\sigma^1$ pour les c\^ones $\sigma \in \EV$ de dimension
$2$. Voici une reformulation plus pr\'ecise du th\'eor\`eme
\ref{theo2}~:

\begin{theobis*}
Pour tout $w \in \Z^2$~:
\[\Vb_w(M) \subset V_{w+(\kappa^1,0)}(M).\]
\end{theobis*}

\subsection{Contr\^ole de la mont\'ee de l'ordre par rapport $V_1$}

Soit $\sigma \in \EV$ un c\^one de dimension maximale et soient $L_1
<L_2$ ses g\'en\'erateurs primitifs. Soit $m\in \Vb_w(M)$ avec $w$
dans $\Z^2$, en particulier $m\in (V^{L_1}_{L_1(w)}(\Dnd)\delta) \cap
(V^{L_2}_{L_2(w)}(\Dnd)\delta) $. Supposons donn\'e $P \in \Dn$ tel
que $P\delta=m$ et $\ord^{L_1}(P)\le L_1(w)$ et tel que $\ord^{L_2}(P)
> L_2(w)$. Alors nous avons montr\'e dans le lemme \ref{lem:cle}
comment construire, en un nombre fini d'\'etapes, un \'el\'ement
$P_\sigma$ tel que $\ord^{L_1}(P_\sigma) \le \ord^{L_1}(P)$
(i.e. l'ordre par rapport \`a $L_1$ n'a pas augment\'e) et
$\ord^{L_2}(P_\sigma) \le L_2(w)$ (i.e. l'ordre par rapport \`a $L_2$
a baiss\'e le plus possible). Nous pouvons nous demander ce qui se
passe pour l'ordre $V_1$ de $P_\sigma$ par rapport \`a celui
de $P$. Nous allons montrer que cet ordre peut augmenter mais de
mani\`ere contr\^ol\'ee. C'est l'objet du lemme suivant~:

\begin{lem}\label{lem:controleV1}
Soit $\sigma$ un c\^one de dimension maximale de $\EV(h(I))$ et $L_1
\ne L_2$ ses g\'en\'erateurs primitifs (il est possible que $L_i
\notin \sigma$). Supposons $V_1 \le L_1 < L_2 \le V_2$.\\
Soient $w\in \Z^2$ et $m \in V^{L_2}_{L_2(w)}(M)$. Soit $P\in \Dnd$
tel que $P\delta=m$ et $\ord^{L_1}(P)\le L_1(w)$ alors on peut
construire $P_\sigma \in \Dnd$ \`a partir de $P$ tel que :
\begin{description}
\item[(i)] $P_\sigma -P \in I$
\item[(ii)] $P_\sigma \in {}^{\sigma}V_w(\Dnd)$, en particulier :
  $\ord^{L_2}(P_\sigma)\le L_2(w)$
\item[(iii)] $\ord^{V_1}(P_\sigma) \le \max\{ \ord^{V_1}(P) \, ,\,
  w_1+\kappa_\sigma^1\}$.
\end{description}
\end{lem}
C'est $(iii)$ qui justifie l'intitul\'e de ce paragraphe.
\begin{proof}[D\'emonstration]
Si $\ord^{L_2}(P)\le L_2(w)$, il suffit de poser $P_\sigma=P$.
On suppose donc que $\ord^{L_2}(P)>L_2(w)$ ce qui entraine
$\ord^{V_1}(\sigma^{L_2}(P)) \le w_1$.

Par hypoth\`ese il existe $P_2 \in \Dnd$ v\'erifiant $P_2\delta=m$ et
$\ord^{L_2}(P_2)\le L_2(w)$. On d\'efinit $H_0=z^{l_0}h(P-P_2)
=z^lh(P)-z^{l_2}h(P_2)$ (il existe des entiers $l_0$, $l$, et $l_2$
satisfaisant une telle \'egalit\'e), $H=z^lh(P)$ et
$H_2=z^{l_2}h(P_2)$. On reprend le d\'ebut de la preuve du lemme
\ref{lem:cle} \`a la diff\'erence qu'on travaille avec la forme $L_2$
au lieu de $L_1$. On consid\`ere donc la division de $H_0$ par la base
standard $Q_1,\ldots,Q_r$ relativement \`a l'ordre $\lhd_{L_2}$ ce qui
donne~:\\
$H_0=\sum_{j=1}^r q_j Q_j$ avec $\ord^{L_2}(H_0) \ge
\ord^{L_2}(q_jQ_j)$. On note $J$ l'ensemble des $j$ dans
$\{1,\ldots,r\}$ pour lesquels l'in\'egalit\'e pr\'ec\'edente est une
\'egalit\'e. On pose alors~:

$W=\sum_{j \in J} \sigma^{L_2}(q_j) Q_j$ et $H'=H-W$.\\
Maintenant, ce qui nous int\'eresse, c'est la diff\'erence entre
$\ord^{V_1}(H)$ et $\ord^{V_1}(H')$. C'est l'objet de ce qui suit~:

\begin{affs*}
\gauche
\begin{description}
\item[(a)] $\ord^{V_1}(W) \le w_1+\kappa_\sigma^1$
\item[(b)] $\ord^{V_1}(W)-\ord^{V_1}(\sigma^{L_2}(W)) \le
  \kappa^1_\sigma$
\end{description}
\end{affs*}
D\'emontrons ces affirmations~:\\
{\bf (a):} Nous avons $\ord^{V_1}(\sigma^{L_2}(W))=
\ord^{V_1}(\sigma^{L_2}(H))=\ord^{V_1}(\sigma^{L_2}(P))$ et ce
dernier est major\'e par $w_1$ donc si {\bf (b)} est vrai il en est de
m\^eme pour {\bf (a)}.\\
{\bf (b):} Comme au paragraphe pr\'ec\'edent, on peut montrer que la
division de $W$ par $\{Q_j,\, j\in J\}$ relativement \`a $\lhd_{L_2}$
donne~: $W=\sum_{j \in J} \sigma^{L_2}(q_j) Q_j$. Par cons\'equent, il
existe $j_1 \in J$ tel que $\exp_{\lhd_{L_2}}(W)=\exp_{\lhd_{L_2}}(m_1
Q_{j_1})$ o\`u $m_1=\mp_{\lhd_{L_2}}(q_{j_1})$ qui est un mon\^one de
$\sigma^{V_1}(\sigma^{L_2}(q_{j_1}))$. En particulier, ceci implique~:
\begin{equation}\label{ident5}
\ord^{V_1}(\sigma^{L_2}(W))= \ord^{V_1}(m_1 \sigma^{L_2}(Q_{j_1})).
\end{equation}
D'autre part,
\[\ord^{V_1}(W) \le \max\{ \ord^{V_1}(\sigma^{L_2}(q_j) Q_j);\, j \in
J \}.\]
Soit alors $j_2 \in J$ tel que $\ord^{V_1}(\sigma^{L_2}(q_{j_2})
Q_{j_2})= \max\{ \ord^{V_1}(\sigma^{L_2}(q_j) Q_j),\,  j\in J\}$. En
prenant $m_2=\mp_{\lhd_{L_2}}(q_{j_2})$, on obtient 
\begin{equation}\label{ident6}
\ord^{V_1}(W)\le \ord^{V_1}(m_2 Q_{j_2})
\end{equation}
Remarquons qu'il est possible qu'on ait $j_1=j_2$. Cependant, on a
toujours~: (voir la figure \ref{fig:affirmations})

\begin{figure}[h!]
\begin{center}
\includegraphics[angle=0, width=11cm]{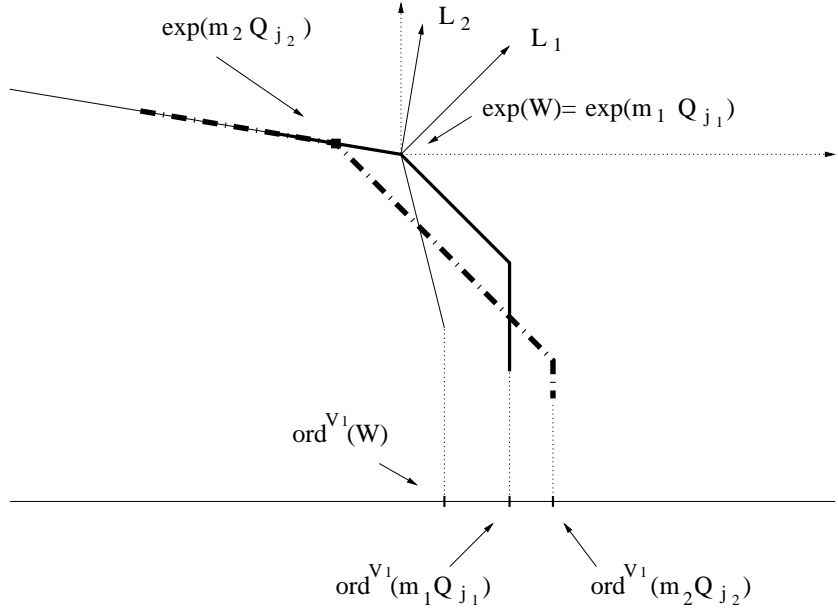}
\caption{Illustration des affirmations}\label{fig:affirmations}
\end{center}
\end{figure}

\begin{aff*}
\gauche
\begin{description}
\item[(c)] $\ord^{V_1}(m_2 \sigma^{L_2}(Q_{j_2})) \le \ord^{V_1}(m_1
  \sigma^{L_2}(Q_{j_1}))$.
\end{description}
\end{aff*}
En utilisant cette affirmation et les identit\'es (\ref{ident5}) et
(\ref{ident6}), nous obtenons~:
\begin{eqnarray*}
\ord^{V_1}(W) & = & \ord^{V_1}(W) - \ord^{V_1}(m_2 Q_{j_2})\\
 & & + \ord^{V_1}(m_2 Q_{j_2}) -\ord^{V_1}(m_2
 \sigma^{L_2}(Q_{j_2}))\\
 & & +\ord^{V_1}(m_2 \sigma^{L_2}(Q_{j_2}))-
 \ord^{V_1}(m_1 \sigma^{L_2}(Q_{j_1}))\\
 & & +\ord^{V_1}(\sigma^{L_2}(W))\\
 & \le & \ord^{V_1}(m_2 Q_{j_2})
 -\ord^{V_1}(m_2 \sigma^{L_2}(Q_{j_2})) +
 \ord^{V_1}(\sigma^{L_2}(W))\\
 & \le & \kappa^1_\sigma+ \ord^{V_1}(\sigma^{L_2}(W))
\end{eqnarray*}
Ceci d\'emontre le point {\bf (b)}. Il ne reste plus qu'\`a
d\'emontrer le point {\bf (c)}.\\

La division de $\sigma^{L_2}(W)$ par $\{\sigma^{L_2}(Q_j),\, j\in J\}$
relativement \`a $\lhd_{L_2}$ donne le r\'esultat suivant~:
$\sigma^{L_2}(W)= \sum_{j\in J} \sigma^{L_2}(q_j) \sigma^{L_2}(Q_j)$. Par
cons\'equent,
\[\exp_{\lhd_{L_2}}\big(\sigma^{L_2}(q_{j_2}) \sigma^{L_2}(Q_{j_2})
\big) \unlhd_{L_2} \exp_{\lhd_{L_2}}(\sigma^{L_2}(W)).\]
Or
\[\exp_{\lhd_{L_2}}(\sigma^{L_2}(W))= \exp_{\lhd_{L_2}}(W)=
\exp_{\lhd_{L_2}}(m_1\sigma^{L_2}(Q_{j_1}))\]
et
\[\exp_{\lhd_{L_2}}\big(\sigma^{L_2}(q_{j_2}) \sigma^{L_2}(Q_{j_2})
\big)= \exp_{\lhd_{L_2}}(m_2\sigma^{L_2}(Q_{j_2}))\]
donc
\[\exp_{\lhd_{L_2}}(m_2\sigma^{L_2}(Q_{j_2})) \unlhd_{L_2}
\exp_{\lhd_{L_2}}(m_1\sigma^{L_2}(Q_{j_1})).\]
Or $\ord^{L_2}(m_2\sigma^{L_2}(Q_{j_2}))=
\ord^{L_2}(m_1\sigma^{L_2}(Q_{j_1}))$ donc
\[\ord^{V_1}(m_2\sigma^{L_2}(Q_{j_2})) \le
\ord^{V_1}(m_1\sigma^{L_2}(Q_{j_1})).\]
Le point {\bf (c)} est d\'emontr\'e.\\
Voyons maintenant comment l'affirmation {\bf (a)} permet de montrer le
troisi\`eme point du lemme. Nous sommes partis de $H$ et nous avons
construit $H'=H-W$. Par {\bf (a)}, nous avons $\ord^{V_1}(H') \le
\max (\ord^{V_1}(H), w_1+\kappa_\sigma^1)$. La suite consiste \`a
faire les m\^emes op\'erations avec $H'$ \`a la place de $H$.
Le dernier \'el\'ement $H_\sigma$ ainsi construit v\'erifie :
$H_\sigma-H \in h(I)$ et $\ord^{V_1}(H_\sigma) \le \max
(\ord^{V_1}(H), w_1+\kappa_\sigma^1)$.\\
On pose alors $P_\sigma={H_\sigma}_{|z=1}$, on a bien $P_\sigma-P \in
I$ et $\ord^{V_1}(P_\sigma) \le \max (\ord^{V_1}(P),
w_1+\kappa_\sigma^1)$. Le lemme est d\'emontr\'e.
\end{proof}

\subsection{Fin de la preuve}

\begin{proof}[D\'emonstration du th\'eor\`eme 2 bis]
Notons $L_0=V_1<\cdots <L_q=V_2$ les \'el\'ements primitifs du
$1$-squelette de $\EV$. Pour chaque $i=1,\ldots,q$, notons
$\sigma_i \in \EV$ le c\^one contenant le c\^one ouvert
engendr\'e par $L_{i-1}$ et $L_i$.\\
Soit $m \in \Vb_w(M)$.\\
Montrons par r\'ecurrence sur $i$ que pour tout $i=0,\ldots,q$, il
existe $T_i \in \Dnd$ v\'erifiant :
\begin{itemize}
\item $T_i \delta=m$
\item $T_i \in V^{L_i}_{L_i(w)}(M)$
\item $\ord^{V_1}(T_i) \le w_1+\kappa^1$
\end{itemize}
Pour $i=0$ : $m\in \Vb_w$ donc en particulier $m \in V^{V_1}_{V_1(w)}$
(notons que $V_1(w)=w_1$) donc il existe $T_0$ tel que $T_0 \delta=m$
et $\ord^{V_1}(T_0)\le w_1 \le w_1+\kappa^1$.\\
Supposons l'assertion vraie au rang $i-1$.\\
On applique le lemme \ref{lem:controleV1} avec $\sigma=\sigma_i$
et $P=T_{i-1}$. On pose alors $T_i=P_\sigma$ (notations du
lemme). D'apr\`es ce m\^eme lemme, $T_i$ v\'erifie :
\begin{itemize}
\item $T_i \delta=m$
\item $T_i \in V^{L_i}_{L_i(w)}(M)$
\item $\ord^{V_1}(T_i) \le \max(\ord^{V_1}(T_{i-1}),w_1+\kappa^1)
  = w_1+\kappa^1$
\end{itemize}
Ainsi, l'assertion est vraie pour tout $i$. En particulier pour $i=q$,
on a : $m= T_q \delta$, $\ord^{V_2}(T_q)\le w_2$ et $\ord^{V_1}(T_q)
\le w_1+\kappa^1$, c'est-\`a-dire $m \in V_{w+(\kappa^1,0)}(M)$.
\end{proof}

\begin{rem}
Le proc\'ed\'e de construction de $(\kappa^1,0)$ montre qu'un tel
$\kappa$ n'est \`a priori pas unique. En effet, en inversant les
r\^oles de $V_1$ et $V_2$, on aurait pu construire un $\kappa$ de la
forme $(0,\kappa^2)$.
\end{rem}
\begin{flushleft}
\textbf{Remerciements}
\nopagebreak

Je remercie mon directeur de th\`ese Michel Granger pour les conseils
donn\'es durant l'\'elaboration de ce travail qui constitue une partie
de ma th\`ese pr\'epar\'ee \`a l'universit\'e d'Angers.
\end{flushleft}

\end{document}